\documentclass[reqno,12pt]{amsart}
\DeclareFontFamily{OML}{script}{}{}
\DeclareFontShape{OML}{script}{m}{it}
{ <5-20> rsfs10 }{}
\DeclareMathAlphabet{\mathscript}{OML}{script}{m}{it}
\usepackage{cite}
\usepackage{color}

\usepackage[pagewise]{lineno}

\usepackage{graphicx,graphics}
\usepackage[colorlinks,breaklinks]{hyperref}
\usepackage{amstext}
\usepackage{pifont}
\usepackage{amsthm}
\usepackage{bm}
\usepackage{amscd}
\usepackage{amsmath}
\usepackage{mathtools}
\usepackage{amssymb}
\usepackage{latexsym}
\usepackage{amsfonts}
\usepackage[notref,final]{showkeys}
\usepackage{esint}
\usepackage{enumerate}
\usepackage{blkarray} 

\usepackage[ansinew]{inputenc}
\usepackage[encapsulated]{CJK}

\ifx\text\undefined
\newcommand{\text}{\mbox}
\fi
\ifx\operatorname\undefined
\newcommand{\operatorname}{\mathop}
\fi

\newcommand{\ki}{{\mbox{\raise.5ex\hbox{$\chi$}}\hspace{.2ex}}}

\renewcommand{\epsilon}{\varepsilon}
\renewcommand{\bar}{\overline}

\newenvironment{eqa}{\begin{equation}%
		\begin{array}{rcl}}{\end{array}\end{equation}}
\newcommand\beqa{\begin{eqa}}
	\newcommand\eeqa{\end{eqa}}

\numberwithin{equation}{section}
\newtheorem{thm}{Theorem}[section]

\newtheorem{rmk}[thm]{Remark}
\newtheorem{prop}[thm]{Proposition}
\newtheorem{defn}[thm]{Definition}

\newcommand{\thmref}[1]{Theorem~\ref{#1}}

\newcommand{\remref}[1]{Remark~\ref{#1}}

\newcommand{\propref}[1]{Proposition~\ref{#1}}

\newcommand{\void}[1]{}

\newcommand{\bR}{{\mathbb R}}  

\begin{document}\begin{CJK}{UTF8}{gkai}

		\title[  Liouville theorem for $\Delta u+f(u)=0$ ]{A New proof of Liouville type theorems for a class of semilinear elliptic equations }
		\author{Chen Guo and Zhengce Zhang}
		\date{\today}
		\address[Chen Guo]{School of Mathematics and Statistics, Xi'an Jiaotong University, Xi'an, 710049, P. R. China}
		\email{jasonchen123@stu.xjtu.edu.cn}
	
		\address[Zhengce Zhang]{School of Mathematics and Statistics, Xi'an Jiaotong University, Xi'an, 710049, P. R. China}
		\email{zhangzc@mail.xjtu.edu.cn}
		\thanks{Corresponding author: Zhengce Zhang}

		\begin{abstract}
			We study  certain typical   semilinear elliptic equations in Euclidean space $\bR^{n}$ or on a closed manifold $M$ with nonnegative Ricci curvature. Our proof is based on a crucial integral identity constructed by the invariant tensor method. Together with  suitable integral estimates, some classical Liouville theorems will be reestablished.
		\end{abstract}
		\maketitle		
		\section{Introduction}\label{Introduction}
		The study of Liouville theorems constitutes a central topic in the theory of partial differential equations. It is well-known that any bounded harmonic function in $\bR^{n}$ is actually constant. For  general Riemannian manifolds, Yau \cite{Y1} proved that every positive harmonic function on a complete Riemannian manifold $M$ with nonnegative Ricci curvature must be constant.
		 Later, Wang and Zhang \cite{WZ} generalized the Liouville theorem to $p$-harmonic functions by deriving the local gradient estimate (also see\cite{KN} under extra assumption on the sectional curvature).
		 Again by gradient estimates, Brighton\cite{B1} establishes the Liouville type result for bounded $f$-harmonic functions on complete smooth metric measure spaces with nonnegative Bakry-Emery curvature.
		 In addition, the moving plane method is also widely used to prove Liouville theorems\cite{CL,CDQ}. There are variations of this approach such as the moving sphere method. Towards the supercritical case, the monotonicity formula is an efficient way to study  Liouville theorems for (stable) solutions (see \cite{DGW,DFP}). To conclude this topic,  we refer to a very recent paper \cite{QS} for a detailed discussion about various approaches to derive Liouville theorems.
		 
		In this paper, we are concerned about the Liouville type theorem for the following semilinear equation 
		\begin{equation}\label{main equ}
			\Delta u+f(u)=0
		\end{equation}
		either in Euclidean space $\bR^{n}$ or on an $n$-dimensional closed manifold $M$ with nonnegative Ricci curvature with $n>2$.
		
		  Equation \eqref{main equ} is fully considered by many researchers and has plenty of conclusions. Owing to the vastness and depth of this topic, we only review some of them here.
		  A typical example of \eqref{main equ} is
		  \begin{equation}\label{power equ}
		  		\Delta u+u^{p}=0.
		  \end{equation}
		  
		    Gidas and Spruck's classical result \cite{GS} shows  that any nonnegative $C^{2}$ solution of \eqref{power equ} in $\bR^{n}$ is identical to zero provided $1< p< \frac{n+2}{n-2}$.
		  It also holds in complete Riemannian manifolds with nonnegative Ricci curvature. In the borderline case $p=\frac{n+2}{n-2}$, Chen and Li \cite{CL} removed the growth restriction at infinity in \cite{GNN} and  proved that \eqref{power equ} has nontrivial radially symmetric solutions  (see also \cite{CGS}). We note that in \cite[Theorem 6.1]{GS}, the authors derive some Liouville type results for more general nonlinearity  $f(x,u)$ with some extra assumptions on the behavior of $f$. 
		  
		  Serrin and Zou \cite{SZ} loosed some technical conditions such as $f(u)>0$ for $u>0$ in \cite{GS}  and extended the conclusion to the quasilinear equation $\Delta_{m}u+f(u)=0$ in $\bR^{n}$. See \cite[Theorem II]{SZ} for details.  They introduced the following subcritical condition. 
		  \begin{defn}\label{subcritical condition1}
		  	The function $f$ is subcritical if $n>m$ and there exists some $\beta\in (0,\frac{nm}{n-m}-1)$ such that
		  	\begin{equation}\label{subcritical formula}
		  		f(u)\ge 0, \quad \beta f(u)-uf'(u)\ge 0, \quad \text{for}\;u>0.
		  	\end{equation}
		  \end{defn}
		  
		  When $m=2$, \eqref{subcritical formula} implies that $u^{-\beta}f(u)$ is non-increasing in $(0,\infty)$ for some $\beta\in (0,\frac{n+2}{n-2})$. Especially, a result involving Laplacian case ($m=2$) is also given.
		  \begin{thm}{\rm(\hspace{-0.03em}\cite[Theorem 4]{SZ})}\label{thm1.1} The following truths hold for \eqref{main equ}.
		  \begin{enumerate}[(i)]
		  \item Let $n=3$ and $f$ is subcritical. Then every solution of \eqref{main equ} is constant.
		  \item Suppose $n\ge 4$. Assume that $f$ is subcritical and the condition 
		  \begin{equation}\label{superpower}
		  	f(u)\ge u^{\alpha}
		  \end{equation}
		  holds for  sufficiently large $u$ with $\alpha>1$. Then \eqref{main equ} only admits the trivial solution $u\equiv 0$. The same conclusion also holds for $\alpha\in (0,1]$ provided $0<\beta\le 1$ in \eqref{subcritical formula}. 
		  \end{enumerate}
		  \end{thm}
		  
		  Afterwards, Li and Zhang \cite{LZ} investigated \eqref{main equ} by the moving sphere method. They found that the solution has two different behaviors, depending on whether $\bar{\lambda}(x)$ is infinite or not (see \cite{LZ}, page 37 for its definition).
		  \begin{thm}{\rm(\hspace{-0.03em}\cite[Theorem 1.3]{	LZ})}\label{thm 1.3} Assume $f$ is locally bounded in $(0,\infty)$ and the function $f(u)u^{-\frac{n+2}{n-2}}$ is non-increasing in $(0,\infty)$. Let $u$ be a solution of \eqref{main equ} in $\bR^{n}$ with $n>2$. Then either $u\equiv a$ for some positive constant $a$  satisfying $f(a)=0$; or for some $b>0$,
		  \begin{equation*}
		  	bu=\left(\frac{\mu}{1+\mu^{2} |x-\bar{x}|^{2}}\right)^{\frac{n-2}{2}},
		  \end{equation*}
		  where $\mu>0,\bar{x}\in \bR^{n}$ and $f(s)s^{-\frac{n+2}{n-2}}\equiv n(n-2)b^{4/(n-2)}$ on $(0,\max_{\bR^{n}} u]$.
		  \end{thm}
		  
		     Since in the critical case $p=\frac{n+2}{n-2}$, no sign assumption of $f$ and also ``super-power'' condition \eqref{superpower}  are required in \thmref{thm 1.3}. One may wonder whether  \thmref{thm1.1} holds by only assuming the subcritical condition for sign-changing $f$.  We attempt to prove them through a
		    totally different way.
		  
		  	In the following discussion, we always suppose that $f\in C^{1}(0,\infty)\bigcap C[0,\infty)$, which may change its sign but has at least one  zero point in $[0,\infty)$. Instead of \eqref{subcritical formula},  we only require $f$ satisfies
		  	\begin{equation}\label{subcritical}
		  		\beta f(u)-uf'(u)\ge 0, \quad \text{for}\;u>0,
		  	\end{equation}
		  	  which implies the function $u^{-\beta} f(u)$ is non-increasing in $(0,\infty)$ for some $\beta\in I$. Here $I$ is an 
		  	  appropriate interval.  For a more intuitive understanding, we conveniently say $f$ is subcritical if $f$ satisfies $\eqref{subcritical}$.
		  	  
		  	  In addition, the solution $u$ of \eqref{main equ} is required at least $C^{3}$. For simplicity, we consider $C^{\infty}$ functions in our theorems. Now we state the main results of this paper. The first theorem investigates the  closed manifold case.
		  	
		  	\begin{thm}\label{mainthm1}
		  		Let $u$ be a smooth positive solution of $\eqref{main equ}$  on an $n$-dimensional closed  manifold $(M,g)$ with nonnegative Ricci curvature. Suppose  $f$ is subcritical for some $\beta\in (0,\frac{n+2}{n-2})$. Then $u$ is constant. The same conclusion also holds if $u^{-\frac{n+2}{n-2}} f(u)$  strictly decreases in $(0,\infty)$.
		  	\end{thm}
		  	
		  	 While considering similar issues in $\bR^{n}$, it turns out that we have to impose the super-power condition about $f$ to handle the integral estimate. It is worthy mentioning that this condition is necessary and can not be removed. 
		  	 
		  	 \begin{thm}\label{mainthm2}
		  	 		Let $u$ be a smooth positive solution of $\eqref{main equ}$ in $\bR^{n}$. Suppose $f$ is subcritical for some $\beta\in (1,\frac{n+2}{n-2})$ and 
		  	 		\begin{equation}\label{condition of mainthm2}
		  	 			|f(u)|\ge c_{1}u^{\beta},
		  	 		\end{equation}
		  	 		where $c_{1}$ is some positive constant dependent on $n,\beta$. Then \eqref{main equ} only admits the trivial solution $u\equiv 0$.
		  	 \end{thm}
		  	
		  	 \begin{thm}\label{mainthm3}
		  	 		Let $u$ be a smooth positive solution of $\eqref{main equ}$ in $\bR^{n}$. Suppose $f$ is subcritical for some $\beta\in (1,\frac{n+2}{n-2})$.   Then the following conclusions hold.
		  	 	\begin{enumerate}[\rm(\roman{enumi})]
		  	 		\item Let $n=3$ and assume $u$ is bounded, then every solution of $u$ is constant.
		  	 		\item Provided $n\ge 4$ and  there exists some $\alpha>1$ such that
		  	 		\begin{equation}\label{condition of mainthm3}
		  	 			f(u)\ge u^{\alpha}  
		  	 		\end{equation}
		  	 		holds for sufficiently large $u$. Then \eqref{main equ} admits no positive solution. Moreover, the only trivial solution is $u\equiv 0$.
		  	 	\end{enumerate}
		  	 \end{thm}
		  	 
		  	 \begin{rmk}
		  	 	\thmref{thm1.1} requires no boundedness assumption of $u$ for $n=3$ with the help of the weak Harnack inequality  and a crucial integral estimate, see \rm(\hspace{-0.03em}\cite[Lemma 3.2]{SZ}) and \rm(\hspace{-0.03em}\cite[Lemma 3.3]{SZ}) respectively. We temporarily lack a similar estimate and 
		  	  	this question will be considered in the future.
		  	 \end{rmk}
		  	 
		  	 \begin{rmk}
		  	    The requirement $\alpha>1$ in Theorem 1.6 {\rm(ii)} originates from the use of Young's inequality.
		  	    
		  	 \end{rmk}
		  	  The soul of our proof is  to derive some important identity.   Gidas and Spruck \cite{GS} made use of the vector field $V$ 
		  	  \begin{equation*}
		  	  	V^{i}=w^{-(n-1)}\left\{w_{ij}w_{j}-\frac{1}{n}w_{i}\Delta w\right\},
		  	  \end{equation*}
		  	  where $w=u^{-\frac{2}{n-2}}$, to  derive some integral identity. It is firstly appeared in  Obata's paper \cite{Obata}.   Different from their approach, we use the following tensor
		  	  \begin{equation*}
		  	  	E_{ij}=u_{ij}+c\frac{u_{i}u_{j}}{u}-\frac{1}{n}\left(\Delta u+c\frac{|\nabla u|^{2}}{u}\right)g_{ij},
		  	  \end{equation*}
		  	  where $c$ is to be determined later. This tensor is firstly introduced by Ma and Wu in \cite{MW} and it is useful to seek for the key identity in Section 2. 
		  	  
		  	  Ma and Ou \cite{MO1} generalized Jerison-Lee's identity \cite{JL1} for the CR Yamabe problem and used their new crucial identity (see Propostion 2.1 of \cite{MO1}) to prove the Liouville theorem for the following equation 
		  	  \begin{equation*}
		  	  	 -\Delta_{\mathbb{H}^{n}}u=2n^{2}u^{q}
		  	  \end{equation*}
		  	  on Heisenberg group $\mathbb{H}^{n}$, where $1<q<\frac{Q+2}{Q-2}$ and $Q=2n+2$ is the homogeneous dimension of $\mathbb{H}^{n}$. Later on, Ma, Ou and Wu \cite{MOW} found the invariant tensor on $\mathbb{H}^{n}$ and sucessfully answered the question asked in \cite{JL1} about the possible theoretical framework to establish such an explicable identity systematically from a perspective. For more relevant studies on the applications of various  (integral) identities, we refer  readers to \cite{CFR,DEL,LM1,LO}  and the references therein.
		  	  
		  	  This paper is organized as follows. In Section 2, we give some necessary notations and establish the significant identity by choosing suitable $c$ in the expression of $E_{ij}$. Then we give the proof of \thmref{mainthm1} in Section 3. Theorems \ref{mainthm2} and \ref{mainthm3} will be proved in Section 4.
\section{ A crucial identity}
\subsection{Notations}
 We introduce some necessary notations at first. Notations in the below adhere Einstein's summation convention, where repeated indices represent sums. Using the standard notation in Riemannian geometry, we  set $\Delta f=f_{i}^{i}$ and $|\nabla u|^{2}=u^{i}u_{i}=g^{ij}u_{i}u_{j}$. $R_{ij}$ denotes  Ricci curvature tensor and $g_{ij}$ represents the common metric tensor. For a second-order tensor $A_{ij}$, its norm $|A_{ij}|$ is defined by $$|A_{ij}|^{2}\triangleq g^{ik}g^{jl}A_{ij}A_{kl}=A_{ij}A^{ij}.$$ From the usual Ricci identity,  the following truth holds
 \begin{equation}\label{Ricci identity}
   f_{ij}^{i}-(\Delta f)_{j}=R_{ij}f^{i}.
 \end{equation}
 
 In addition, $C(\cdot)$ always denotes some positive constant depending on parameters in the bracket, which may change the value from line to line.
 
 \subsection{Derivation of the identity}
 Our crucial identity is given as follows.
 \begin{prop}\label{prop1}
 	\begin{equation}\label{keyidentity}
 			\begin{aligned}
 			u^{\frac{2\beta}{n+2}}\left(u^{-\frac{2\beta}{n+2}}E_{ij}u^{j}\right)^{i}=&|E_{ij}|^{2}+\frac{n-1}{n+2}\beta\left(1-\frac{n-2}{n+2}\beta\right)\frac{|\nabla u|^{4}}{u^{2}}\\
 			&+R_{ij}u^{i}u^{j}-\frac{n-1}{n}\left(f'(u)-\beta\frac{f(u)}{u}\right)|\nabla u|^{2}.
 		\end{aligned}
 	\end{equation}
 \end{prop}
 \noindent{\bf Proof of Proposition 2.1}\\
   Recall
   \begin{equation}\label{tensorEij}
   	 E_{ij}=u_{ij}+c\frac{u_{i}u_{j}}{u}-\frac{1}{n}\left(\Delta u+c\frac{|\nabla u|^{2}}{u}\right)g_{ij}.
   \end{equation}
  
   Then a straightforward calculation shows
   \begin{equation}\label{div1}
   	E_{ij}^{i}=\frac{n-1}{n}(\Delta u)_{j}+R_{ij}u^{i}+c\frac{\Delta u u_{j}}{u}+\frac{n-2}{n}c\frac{u_{ij}u^{i}}{u}-\frac{n-1}{n}c\frac{
   	|\nabla u|^{2}}{u^{2}}u_{j},
   \end{equation}
   where we have used \eqref{Ricci identity}.
   
    By using the fact $(\Delta u)_{j}=-f'(u)u_{j}$ and replacing $u_{ij}$ by $E_{ij}$,  a combination of \eqref{tensorEij} and \eqref{div1} leads to
    \begin{equation}\label{div2}
    	\begin{aligned}
    		E_{ij}^{i}=&-\frac{n-1}{n}f'(u)u_{j}+R_{ij}u^{i}+\frac{(n-1)(n+2)}{n^{2}}c\frac{\Delta u u_{j}}{u}\\
    		&+\frac{n-2}{n}c\frac{E_{ij}u^{i}}{u}-\frac{n-1}{n}c\left(1+\frac{n-2}{n}c\right)\frac	{|\nabla u|^{2}}{u^{2}}u_{j}.
    	\end{aligned}
   \end{equation}
   Interpolate the term $$\frac{n-1}{n}\beta\frac{f(u)}{u}u_{j}$$ into \eqref{div2}, it yields
   \begin{equation}\label{div3}
   	\begin{aligned}
   		E_{ij}^{i}=&-\frac{n-1}{n}\left(f'(u)-\beta\frac{f(u)}{u}\right)u_{j}+R_{ij}u^{i}+\frac{n-2}{n}c\frac{E_{ij}u^{i}}{u}\\
   		&+\left(\frac{(n-1)(n+2)}{n^{2}}c+\frac{n-1}{n}\beta\right)\frac{\Delta u u_{j}}{u}-\frac{n-1}{n}c\left(1+\frac{n-2}{n}c\right)\frac	{|\nabla u|^{2}}{u^{2}}u_{j}.
   	\end{aligned}
   \end{equation}
   In order to eliminate the mixed term $\frac{\Delta u u_{j}}{u}$, we set $	c=-\frac{n\beta}{n+2}$ and \eqref{div3} becomes
   \begin{equation}\label{divfinal}
   	\begin{aligned}
   		E_{ij}^{i}=&-\frac{n-1}{n}\left(f'(u)-\beta\frac{f(u)}{u}\right)u_{j}+R_{ij}u^{i}\\&-\frac{n-2}{n+2}\beta\frac{E_{ij}u^{i}}{u}+\frac{n-1}{n+2}\beta(1-\frac{n+2}{n-2}\beta)\frac{|\nabla u|^{2}}{u^{2}}u_{j}.
   	\end{aligned}
   \end{equation}
   Due to trace free property of $E_{ij}$, we combine $\eqref{tensorEij}$ and \eqref{divfinal} to obtain
    \begin{equation}\label{equ 2.8}
    	\begin{aligned}
    	\left(E_{ij}u^{j}\right)^{i}=&|E_{ij}|^{2}+\frac{2\beta}{n+2}\frac{E_{ij}u^{i}u^{j}}{u}-\frac{n-1}{n}\left(f'(u)-\beta\frac{f(u)}{u}\right)|\nabla u|^{2}\\
    	&+\frac{n-1}{n+2}\beta(1-\frac{n-2}{n+2}\beta)\frac{|\nabla u|^{4}}{u^{2}}+R_{ij}u^{i}u^{j}
    	\end{aligned}
   \end{equation}

   With the purpose of removing the second term in the right hand side of \eqref{equ 2.8}, we introduce a certain power of $u$ to accomplish it. Notice that
    \begin{equation}\label{equ 2.9}
    	u^{\frac{2\beta}{n+2}}\left(u^{-\frac{2\beta}{n+2}}E_{ij}u^{j}\right)^{i}=\left(E_{ij}u^{j}\right)^{i}-\frac{2\beta}{n+2}\frac{E_{ij}u^{i}u^{j}}{u}.
    \end{equation}
   We immediately get \eqref{keyidentity} by substituting \eqref{equ 2.8} into \eqref{equ 2.9}. \qed 
   

\section{Proof of Theorem 1.4}
   We  prove \thmref{mainthm1} at first. It follows from \propref{prop1} that
  	\begin{equation}\label{identity}
  	\begin{aligned}
  		\left(u^{-\frac{2\beta}{n+2}}E_{ij}u^{j}\right)^{i}=&u^{-\frac{2\beta}{n+2}}|E_{ij}|^{2}+\frac{n-1}{n+2}\beta\left(1-\frac{n-2}{n+2}\beta\right)u^{-\frac{2\beta}{n+2}-2}|\nabla u|^{4}\\
  		&+u^{-\frac{2\beta}{n+2}}R_{ij}u^{i}u^{j}-\frac{n-1}{n}\left(f'(u)-\beta\frac{f(u)}{u}\right)u^{-\frac{2\beta}{n+2}}|\nabla u|^{2}.
  	\end{aligned}
  \end{equation}
  
  In the treatment of the closed manifold, we just integrate \eqref{identity} over $M$ and use the divergence theorem to obtain
  \begin{equation}\label{final1}
      \begin{aligned}
      		0=&\int_{M}\Big[u^{-\frac{2\beta}{n+2}}|E_{ij}|^{2}+\frac{n-1}{n+2}\beta\left(1-\frac{n-2}{n+2}\beta\right)u^{-\frac{2\beta}{n+2}-2}|\nabla u|^{4}\\
      		&+u^{-\frac{2\beta}{n+2}}R_{ij}u^{i}u^{j}-\frac{n-1}{n}\left(f'(u)-\beta\frac{f(u)}{u}\right)u^{-\frac{2\beta}{n+2}}|\nabla u|^{2} \Big].
      \end{aligned}
   \end{equation}
   
   Notice that the coefficient $\beta\left(1-\frac{n-2}{n+2}\beta\right)$ keeps positive when $\beta\in (0,\frac{n+2}{n-2})$. Using \eqref{subcritical}  and the non-negativity of the Ricci curvature, \eqref{final1} implies
   \begin{equation}
   	0\ge \int_{M} \frac{n-1}{n+2}\beta\left(1-\frac{n-2}{n+2}\beta\right)u^{-\frac{2\beta}{n+2}-2}|\nabla u|^{4},
   \end{equation}
   which enforces $u$ to be constant.
   
   Furthermore, if $u^{-\frac{n+2}{n-2}}f(u)$ is strictly decreasing in $(0,\infty)$, we can easily get
   \begin{equation}\label{equ 3.4}
   	f'(u)-\frac{n+2}{n-2}\frac{f(u)}{u}<0.
   \end{equation}
   
   Setting $\beta=\frac{n+2}{n-2}$ in \eqref{final1}, the second term in the right hand side vanishes and we retain that
   \begin{equation}\label{final2}
   	0\ge -\int_{M} \frac{n-1}{n}\left(f'(u)-\frac{n+2}{n-2}\frac{f(u)}{u}\right)u^{-\frac{2}{n-2}}|\nabla u|^{2} .
   \end{equation}
   As a consequence of \eqref{equ 3.4} and \eqref{final2}, $u$ should be constant.  This ends the proof.
   \section{The Euclidean case: Proofs of Theorems 1.5 and  1.6}
   Next, we deal with the Euclidean case. At this point, the covariant derivative degenerates into the usual partial derivative in $\bR^{n}$, $g_{ij}=\delta_{ij}$ and $R_{ij}=0$. We always assume  $1< \beta<\frac{n+2}{n-2}$ in this section. Since we are going to derive some integral estimates, in order to be consistent with the standard notation in partial differential equations, we will uniformly adopt subscripts below.
   \subsection{Some integral estimates} We prove the following integral estimates at the beginning. They are necessary in our proof. 
   At first, we choose a cut-off function $\eta\in C^{\infty}_{0}(R^{n})$ satisfying 
  \begin{equation*}
  	 0\le\eta\le 1,\; \eta\equiv\,1\;\text{in}\;B_{R},\; \eta=0 \;\text{outside}\; B_{2R}\;\,\text{and}\;|\nabla \eta|\le \frac{C(n)}{R}.
  \end{equation*}
  Multiplying both sides of \eqref{identity} by $\eta^{\gamma}$ and integrating over $\bR^{n}$, we obtain
  \begin{equation}\label{indentity for RN}
  	\begin{aligned}
  &\int_{\bR^{n}}u^{-\frac{2\beta}{n+2}}|E_{ij}|^{2}\eta^{\gamma}+\int_{\bR^{n}}\frac{n-1}{n+2}\beta\left(1-\frac{n-2}{n+2}\right)u^{-\frac{2\beta}{n+2}-2}|\nabla u|^{4}\eta^{\gamma}\\
  =&\int_{\bR^{n}}\left(u^{-\frac{2\beta}{n+2}}E_{ij}u_{j}\right)_{i}\eta^{\gamma}+\int_{\bR^{n}}\frac{n-1}{n}\left(f'(u)-\beta\frac{f(u)}{u}\right)u^{-\frac{2\beta}{n+2}}|\nabla u|^{2}\eta^{\gamma}.
  	\end{aligned}
  \end{equation}
  
  Here $\gamma$ relies on $n,\beta$ and will be chosen later. We denote $A=\beta\left(1-\frac{n-2}{n+2}\right)>0$. For simplicity, in the below the integration region is always ignored when integrating in the whole space. A
  direct integration by parts yields that
  \begin{equation*}
  	\int\left(u^{-\frac{2\beta}{n+2}}E_{ij}u_{j}\right)_{i}\eta^{\gamma}
  =-\gamma \int u^{-\frac{2\beta}{n+2}}E_{ij}u_{j}\eta^{\gamma-1}\eta_{i}.
  \end{equation*}
  Using  Cauchy-Schwarz's inequality, we have
  \begin{equation}\label{equ 4.2}
  	\begin{aligned}
  		&-\gamma \int u^{-\frac{2\beta}{n+2}}E_{ij}u_{j}\eta^{\gamma-1}\eta_{i}\\
  		&\le
  		 \frac{1}{2}\int u^{-\frac{2\beta}{n+2}}|E_{ij}|^{2}\eta^{\gamma}
  		+C(n,\beta)\int u^{-\frac{2\beta}{n+2}} |\nabla u|^{2}\eta^{\gamma-2}|\nabla \eta|^{2}\\
  		&\le  \frac{1}{2}\int u^{-\frac{2\beta}{n+2}}|E_{ij}|^{2}\eta^{\gamma}
  		+\frac{A}{2}\int u^{-\frac{2\beta}{n+2}}|\nabla u|^{4}\eta^{\gamma}
  		+C(n,\beta)\int u^{-\frac{2\beta}{n+2}+2}\eta^{\gamma-4}|\nabla \eta|^{4}.
  	\end{aligned}
  \end{equation}
  With the help of \eqref{subcritical}, we combine \eqref{indentity for RN} and \eqref{equ 4.2} to get
  \begin{equation}\label{first integral inequ}
  	\int u^{-\frac{2\beta}{n+2}-2}|\nabla u|^{4}\eta^{\gamma}\le C(n,\beta)\int u^{-\frac{2\beta}{n+2}+2}\eta^{\gamma-4}|\nabla \eta|^{4}.
  \end{equation}
  From equation \eqref{main equ}, we have
  \begin{equation}\label{identity from mainequ}
            \int u^{-\frac{2\beta}{n+2}} f^{2}(u)\eta^{\gamma}=-\int u^{-\frac{2\beta}{n+2}} f(u)\eta^{\gamma}\Delta u.
  \end{equation}
  Integrate by parts to obtain
  \begin{equation}\label{integration by parts1}
  	\begin{aligned}
  		 -\int u^{-\frac{2\beta}{n+2}} f(u)\eta^{\gamma}\Delta u
  		 =&\int \gamma\eta^{\gamma-1}u^{-\frac{2\beta}{n+2}}f(u)\nabla\eta\cdot \nabla u+\int u^{-\frac{2\beta}{n+2}} f'(u)|\nabla u|^{2}\eta^{\gamma}\\
  		 &-\int \frac{2\beta}{n+2}u^{-\frac{2\beta}{n+2}-1} f(u)|\nabla u|^{2}\eta^{\gamma}.
  	\end{aligned}
  \end{equation}
 By reapplying \eqref{subcritical} and combining it with \eqref{identity from mainequ} and \eqref{integration by parts1}, we derive that
 \begin{equation}\label{equ 4.6}
 		 \int u^{-\frac{2\beta}{n+2}} f^{2}(u)\eta^{\gamma}\le \frac{n\beta}{n+2}\int u^{-\frac{2\beta}{n+2}-1} f(u)|\nabla u|^{2}\eta^{\gamma}+\gamma\int \eta^{\gamma-1}u^{-\frac{2\beta}{n+2}}f(u)\nabla\eta\cdot \nabla u.
 \end{equation}
 Note that
 \begin{equation}\label{equ 4.7}
 	\begin{aligned}
 		\left| \frac{n\beta}{n+2}\int u^{-\frac{2\beta}{n+2}-1} f(u)|\nabla u|^{2}\eta^{\gamma} \right|&\le \int u^{-\frac{2\beta}{n+2}-1} |f(u)| |\nabla u|^{2}\eta^{\gamma}\\
 		&\le \frac{1}{4}  \int u^{-\frac{2\beta}{n+2}} f^{2}(u)\eta^{\gamma}
 		+C(n,\beta)\int  u^{-\frac{2\beta}{n+2}+2}\eta^{\gamma-4}|\nabla \eta|^{4},
 	\end{aligned}
 \end{equation}
 \begin{equation}\label{equ 4.8}
 	\begin{aligned}
 		\left| \gamma\int \eta^{\gamma-1}u^{-\frac{2\beta}{n+2}}f(u)\nabla\eta\cdot \nabla u \right|&\le  \gamma\int u^{-\frac{2\beta}{n+2}}|f(u)| |\nabla\eta| |\nabla u| \eta^{\gamma-1}\\
 		&\le \frac{1}{4}  \int u^{-\frac{2\beta}{n+2}} f^{2}(u)\eta^{\gamma}
 		+\gamma^{2}\int  u^{-\frac{2\beta}{n+2}}|\nabla u|^{2}|\nabla \eta|^{2}\eta^{\gamma-2},
 	\end{aligned}
 \end{equation}
 and
 \begin{equation}\label{equ 4.9}
 	\begin{aligned}
 		\gamma^{2}\int  u^{-\frac{2\beta}{n+2}}|\nabla u|^{2}|\nabla \eta|^{2}\eta^{\gamma-2}&\le \frac{1}{4} \int u^{-\frac{2\beta}{n+2}-2}|\nabla u|^{4}\eta^{\gamma}+\gamma^{4} \int  u^{-\frac{2\beta}{n+2}+2}\eta^{\gamma-4}|\nabla \eta|^{4}\\
 		&\le C(n,\beta)\int  u^{-\frac{2\beta}{n+2}+2}\eta^{\gamma-4}|\nabla \eta|^{4},
 	\end{aligned}
 \end{equation}
  where we have used \eqref{first integral inequ} in the second inequality of \eqref{equ 4.9}.
  
   \noindent Substituting \eqref{equ 4.7}, \eqref{equ 4.8} and \eqref{equ 4.9} into \eqref{equ 4.6}, it yields
   \begin{equation}\label{second integral inequality}
   	\int u^{-\frac{2\beta}{n+2}} f^{2}(u)\eta^{\gamma}\le C(n,\beta)\int  u^{-\frac{2\beta}{n+2}+2}\eta^{\gamma-4}|\nabla \eta|^{4}.
   \end{equation}
   
   We have now rigorously established the target integral inequalites via the preceding analysis.
   \begin{rmk}\label{rmk4.1}
   	Both \eqref{first integral inequ} and \eqref{second integral inequality} still hold for a complete noncompact manifold with nonnegative Ricci curvature. Instead of \eqref{indentity for RN}, we have
   	\begin{equation}\label{identity for mfd}
   		\begin{aligned}
   			\int_{M}u^{-\frac{2\beta}{n+2}}R_{ij}u^{i}u^{j}&+\int_{M}u^{-\frac{2\beta}{n+2}}|E_{ij}|^{2}\eta^{\gamma}+\int_{M}\frac{n-1}{n+2}\beta\left(1-\frac{n-2}{n+2}\right)u^{-\frac{2\beta}{n+2}-2}|\nabla u|^{4}\eta^{\gamma}\\
   			&=\int_{M}\left(u^{-\frac{2\beta}{n+2}}E_{ij}u_{j}\right)_{i}\eta^{\gamma}+\int_{M}\frac{n-1}{n}\left(f'(u)-\beta\frac{f(u)}{u}\right)u^{-\frac{2\beta}{n+2}}|\nabla u|^{2}\eta^{\gamma}.
   		\end{aligned}
   	\end{equation}
   	
   	The first term in the left hand side of \eqref{identity for mfd} could be dropped due to the condition $Ric(M)\ge 0$. Then a comparable process also produces these two formulae.
   \end{rmk}
   \begin{rmk}
   	It can be easily seen that
   	\begin{equation*}
   			\int u^{-\frac{2\beta}{n+2}-2}|\nabla u|^{4}\eta^{\gamma}+\int u^{-\frac{2\beta}{n+2}} f^{2}(u)\eta^{\gamma}\le C(n,\beta)\int  u^{-\frac{2\beta}{n+2}+2}\eta^{\gamma-4}|\nabla \eta|^{4}.
   	\end{equation*}
   	 The above formula is a special case of \rm(\hspace{-0.03em}\cite[Lemma 5.3]{he2025optimalliouvilletheoremslaneemden}).
   \end{rmk}
   \subsection{Proof of Theorem 1.5} From \eqref{second integral inequality} and \eqref{condition of mainthm2}, we immediately get
   \begin{equation}\label{integrate estimate for thm1.4}
   			\int u^{\frac{2n+2}{n+2}\beta} \eta^{\gamma}\le C(n,\beta)\int  u^{-\frac{2\beta}{n+2}+2}\eta^{\gamma-4}|\nabla \eta|^{4}.
   \end{equation}
   
   Now the same argument as in \cite{MW} is available.  We sketch the proof here. Since  $1<\beta<\frac{n+2}{n-2}$ suggests $2-\frac{2\beta}{n+2}>2-\frac{2}{n-2}>0$ and $\frac{2n+2}{n+2}\beta+\frac{2\beta}{n-2}-2=2\beta-2>0$, Young's inequality implies that
   \begin{equation*}
   	     C(n,\beta)\int  u^{-\frac{2\beta}{n+2}+2}\eta^{\gamma-4}|\nabla \eta|^{4}\le \frac{1}{2} \int u^{\frac{2n+2}{n+2}\beta} \eta^{\gamma}
   	     +C_{1}(n,\beta)\int \eta^{\gamma-\frac{4(n+1)\beta}{(n+2)(\beta-1)}}
   	     |\nabla \eta|^{\frac{4(n+1)\beta}{(n+2)(\beta-1)}}.
   \end{equation*}
   Consequently, we choose $\gamma=1+\frac{4(n+1)\beta}{(n+2)(\beta-1)}$ to obtain that
   \begin{equation}\label{finalestimate1}
   		\int_{B_{R}} u^{\frac{2n+2}{n+2}\beta} \le C_{1}(n,\beta)R^{n-\frac{4(n+1)\beta}{(n+2)(\beta-1)}}.
   \end{equation} 
   Since
   \begin{equation*}
   	n-\frac{4(n+1)\beta}{(n+2)(\beta-1)}=\frac{(n^{2}-2n-4)\beta-n(n+2)}{(n+2)(\beta-1)}
   \end{equation*} 
   and
   \begin{equation*}
   	(n^{2}-2n-4)\beta-n(n+2)<-\frac{4(n+2)}{n-2}<0,
   \end{equation*} we know the right hand side of \eqref{finalestimate1} tends to zero as $R\to \infty$.
   This contracts the positivity of $u$ and suggests \eqref{main equ} admits no positive solution. Recall $f(u)$ possesses at least one root in $[0,\infty)$, \eqref{condition of mainthm2} means the only zero point of $f$ is $u=0$. Thus the only trivial solution of \eqref{main equ} is $u\equiv 0$. This completes the proof. \qed
   
   \subsection{Proof of Theorem 1.6}
   	  We first prove (ii). Since $f$ is subcritical, \eqref{subcritical} implies 
   	   \begin{equation}\label{subcritical2}
   	   	\left(u^{-\beta}f(u)\right)'\le 0,
   	   \end{equation}
   	   where the derivative is with respect to $u$. Meanwhile, \eqref{condition of mainthm3} suggests there exist some  $u_{0}$ sufficiently large such that 	$f(u)\ge u^{\alpha}$ holds for $u\ge u_{0}$ (Without loss of generality, we choose $u_{0}>1$). Therefore $f(u_{0})>0$. Integrate \eqref{subcritical2} from $u_{0}$ to $u$ to get
   	   \begin{equation}\label{inequ for u large}
   	   	    u^{\alpha}\le f(u)\le \frac{f(u_{0})}{u_{0}^{\beta}}u^{\beta},\qquad u\ge u_{0}.
   	   \end{equation}
   	   
   	   It is worthy mentioning that \eqref{inequ for u large} suggests $\alpha\le \beta$ (See \cite{SZ}, page 102). Proceeding in a similar manner, we also derive 
   	   \begin{equation}\label{inequ for u small}
   	    f(u)\ge \frac{f(u_{0})}{u_{0}^{\beta}}u^{\beta},\qquad u< u_{0}.
   	   \end{equation}
   	   
   	    We use a trick in the proof of \cite[Theorem 2.8]{he2025optimalliouvilletheoremslaneemden}. Decompose $\bR^{n}$ into two disjoint parts $\{u<u_{0}\}$ and $\{u\ge u_{0}\}$. More precisely,
   	    \begin{equation*}
   	    	\{u\ge u_{0}\}\triangleq \{x\in \bR^{n}:\; u(x)\ge u_{0}\}.
   	    \end{equation*}
   	    Now we aim to derive an integral estimate similar to \eqref{second integral inequality}. Note that
   	    $$ -\frac{2\beta}{n+2}+2=\frac{2\beta}{n+2}\cdot \frac{n+2-\beta}{\beta}$$
   	    and 
   	    \begin{equation*}
   	    	p_{1}=\frac{(n+1)\beta}{n+2-\beta}>1,\qquad q_{1}=\frac{(n+1)\beta}{(n+2)(\beta-1)}>0.
   	    \end{equation*}
   	    By virtue of H\"{o}lder's inequality with the respective exponent pair $(p_{1},q_{1})$, it yields 
   	    \begin{equation}\label{estimate for small}
   	    	\begin{aligned}
   	    		\int_{\{u<u_{0}\}} u^{-\frac{2\beta}{n+2}+2} \eta^{\gamma-4} |\nabla\eta|^{4}
   	    		&\le \left( \int_{\{u<u_{0}\}} u^{\frac{(2n+2)\beta}{n+2}} \eta^{\gamma} \right)^{\frac{n+2-\beta}{(n+1)\beta}} \\
   	    		&\quad \times \left( \int_{\{u<u_{0}\}} \eta^{\gamma-\frac{(4n+4)\beta}{(n+2)(\beta-1)}} |\nabla \eta|^{\frac{(4n+4)\beta}{(n+2)(\beta-1)}} \right)^{\frac{(n+2)(\beta-1)}{(n+1)\beta}} \\
   	    		&\le C(n,\beta) \left( \int_{\{u<u_{0}\}} f^{2}(u) u^{-\frac{2\beta}{n+2}} \eta^{\gamma} \right)^{\frac{n+2-\beta}{(n+1)\beta}} \\
   	    		&\quad \times \left( \int_{\{u<u_{0}\}} \eta^{\gamma-\frac{(4n+4)\beta}{(n+2)(\beta-1)}} |\nabla \eta|^{\frac{(4n+4)\beta}{(n+2)(\beta-1)}} \right)^{\frac{(n+2)(\beta-1)}{(n+1)\beta}}. \\
   	    	\end{aligned}
   	    \end{equation}
   	  We analogously apply H\"{o}lider's inequality with the exponents
   	  $$ p_{2}=\frac{\alpha-\frac{\beta}{n+2}}{1-\frac{\beta}{n+2}}>1 \quad
   	  \text{and}\quad q_{2}=\frac{\alpha-\frac{\beta}{n+2}}{\alpha-1}>1 $$
   	  and derive the following estimate
   	  \begin{equation}\label{estimate for large}
   	  	\begin{aligned}
   	  		\int_{\{u\ge u_{0}\}} u^{-\frac{2\beta}{n+2}+2} \eta^{\gamma-4} |\nabla\eta|^{4}
   	  		&\le \left( \int_{\{u\ge u_{0}\}} u^{2\alpha-\frac{2\beta}{n+2}} \eta^{\gamma} \right)^{\frac{n+2-\beta}{(n+2)\alpha-\beta}} \\
   	  		&\quad \times \left( \int_{\{u\ge u_{0}\}} \eta^{\gamma-4\frac{(n+2)\alpha-\beta}{(n+2)(\alpha-1)}} |\nabla \eta|^{4\frac{(n+2)\alpha-\beta}{(n+2)(\alpha-1)}} \right)^{\frac{(n+2)(\alpha-1)}{(n+2)\alpha-\beta}} \\
   	  	&\le \left( \int_{\{u\ge u_{0}\}} f^{2}(u)u^{-\frac{2\beta}{n+2}} \eta^{\gamma} \right)^{\frac{n+2-\beta}{(n+2)\alpha-\beta}} \\
   	  	&\quad \times \left( \int_{\{u\ge u_{0}\}} \eta^{\gamma-4\frac{(n+2)\alpha-\beta}{(n+2)(\alpha-1)}} |\nabla \eta|^{4\frac{(n+2)\alpha-\beta}{(n+2)(\alpha-1)}} \right)^{\frac{(n+2)(\alpha-1)}{(n+2)\alpha-\beta}}.
   	  	\end{aligned}
   	  \end{equation}
   A combination of \eqref{second integral inequality}, \eqref{estimate for small} and \eqref{estimate for large} leads to
	\begin{equation}\label{equ 4.19}
		\begin{aligned}
				\int u^{-\frac{2\beta}{n+2}} f^{2}(u)\eta^{\gamma}
				&\le C(n,\beta) \left( \int_{\{u<u_{0}\}} f^{2}(u) u^{-\frac{2\beta}{n+2}} \eta^{\gamma} \right)^{\frac{n+2-\beta}{(n+1)\beta}} \\
				&\quad \times \left( \int_{\{u<u_{0}\}} \eta^{\gamma-\frac{(4n+4)\beta}{(n+2)(\beta-1)}} |\nabla \eta|^{\frac{(4n+4)\beta}{(n+2)(\beta-1)}} \right)^{\frac{(n+2)(\beta-1)}{(n+1)\beta}} \\
				 &\quad +C(n,\beta) \left( \int_{\{u\ge u_{0}\}} f^{2}(u)u^{-\frac{2\beta}{n+2}} \eta^{\gamma} \right)^{\frac{n+2-\beta}{(n+2)\alpha-\beta}} \\
				 &\quad \times \left( \int_{\{u\ge u_{0}\}} \eta^{\gamma-4\frac{(n+2)\alpha-\beta}{(n+2)(\alpha-1)}} |\nabla \eta|^{4\frac{(n+2)\alpha-\beta}{(n+2)(\alpha-1)}} \right)^{\frac{(n+2)(\alpha-1)}{(n+2)\alpha-\beta}}.
		\end{aligned}
	\end{equation}
  	We choose $$\gamma=\max\left\{ 1+\frac{(4n+4)\beta}{(n+2)(\beta-1)},1+4\frac{(n+2)\alpha-\beta}{(n+2)(\alpha-1)} \right\}  $$
  	in \eqref{equ 4.19} and use Young's inequality. By the definition of $\eta$, it provides
  	\begin{equation}\label{finalestimate2}
  		\int_{B_{R}}  u^{-\frac{2\beta}{n+2}} f^{2}(u) \le C(n,\alpha,\beta)\left[R^{n-4\frac{(n+1)\beta}{(n+2)(\beta-1)}}+R^{n-4\frac{(n+2)\alpha-\beta}{(n+2)(\alpha-1)}}\right].
  	\end{equation}
  	
  We claim that the exponents in the right hand side of \eqref{finalestimate2} are both negative. For $n=3$, the first exponent becomes
  \begin{equation*}
  	3-\frac{16\beta}{5(\beta-1)}=\frac{-\beta-15}{5(\beta-1)}<0.
  \end{equation*}
  For $n\ge 4$, we have
  \begin{equation*}
  	n-4\frac{(n+1)\beta}{(n+2)(\beta-1)}< n-\frac{4(n+1)}{n+2}\left(1+\frac{n-2}{4}\right)=-1<0.
  \end{equation*}
  Thus the first exponent is negative. Concerning the second index, we have
  \begin{equation*}
  	n-4\frac{(n+2)\alpha-\beta}{(n+2)(\alpha-1)}=\frac{n(n+2)(\alpha-1)-4(n+2)\alpha+4\beta}{(n+2)(\alpha-1)}.
  \end{equation*}
  A straightforward computation shows
  \begin{equation}\label{determiant}
  	n(n+2)(\alpha-1)-4(n+2)\alpha+4\beta=(n^{2}-2n-8)\alpha-n^{2}-2n+4\beta.
  \end{equation}
  For $n=3$, we use the fact $\beta<5$ and $\alpha>1$ and \eqref{determiant} becomes
  \begin{equation*}
  	 -5\alpha-9-6+4\beta<-5\alpha+5<0.
  \end{equation*}
  For  $n\ge 4$, since $n^{2}-2n-8\ge 0$, we have that
  \begin{align*}
  	(n^{2}-2n-8)\alpha-n^{2}-2n+4\beta&\le (n^{2}-2n-4)\beta-n^{2}-2n \\
  	&<  (n^{2}-2n-4)\frac{n+2}{n-2}-n^{2}-2n\\
  	&<0.
  \end{align*}  
  
   Therefore, the negativity of \eqref{determiant} implies that the second exponent is also negative. Letting $R\to \infty$, \eqref{finalestimate2}
   suggests that $f(u)\equiv 0$ which contracts \eqref{condition of mainthm3}.
   Consequently, \eqref{main equ} admits no positive solution. Furthermore, \eqref{inequ for u large} suggests $f$ has no zero point in $[u_{0},\infty)$. If $f(u)$  vanishes at a nonzero  $\xi\in (0,u_0)$,
    then take $u=\xi$ in \eqref{inequ for u small} to derive
    $$  0=f(\xi)\ge \frac{f(u_{0})}{u_{0}^{\beta}}\xi^{\beta}>0,$$
   which is impossible. The above discussion indicates that the unique zero point of $f$ is $u=0$ and the only trivial solution of \eqref{main equ} is $u\equiv 0$. This proves (ii).
  
   For the special case $n=3$, a stronger conclusion holds if  $u$ is bounded, i.e. $u\le L$. We note that the derivation of \eqref{second integral inequality} does not require super-power condition like \eqref{superpower}. In this case, directly from \eqref{second integral inequality},
   we set $\gamma=4$, $n=3$  and easily get 
   \begin{equation}\label{finalestimate3}
   	   \int_{B_{R}} u^{-\frac{2\beta}{n+2}} f^{2}(u)\le C(n,\beta,L)R^{-1}.
   \end{equation}
   Hence $f(u)\equiv 0$ and the classical Liouville theorem for harmonic functions implies that $u$ is constant.  This proves (i). \qed
   \begin{rmk}
   	   The arguments presented in this section can be extended to general complete noncompact Riemannian manifolds with nonnegative Ricci curvature with a minor adjustment by \remref{rmk4.1} and the standard volume comparison theorem.
   \end{rmk}
   
	 \noindent{\bf Acknowledgments.}
		The authors would thank  Professor Xi-Nan Ma for introducing us this interesting tool. This work was partially supported by the National Natural Science Foundation of China (No. 12271423) and the Shaanxi Fundamental Science Research Project for Mathematics and Physics (No. 23JSY026).

	\bibliographystyle{amsplain}
	\bibliography{reference}
\end{CJK}
\end{document}